\documentclass[twocolumn]{autart}    

\usepackage[T1]{fontenc}
\usepackage{palatino} 
\usepackage{mathrsfs}
\usepackage{graphicx}
\usepackage{mathtools}
\usepackage{amsmath} 
\usepackage{amssymb}  
\usepackage{color}
\usepackage{enumerate}
\usepackage{multirow}
\usepackage{nomencl}
\usepackage{array,url}

\newcommand{\bmat}[1]{\begin{bmatrix}#1\end{bmatrix}}
\usepackage{chngcntr}
\usepackage{apptools}
\AtAppendix{\counterwithin{lemma}{section}}
\AtAppendix{\counterwithin{assumption}{section}}
\AtAppendix{\counterwithin{proposition}{section}}

\newtheorem{lemma}{Lemma}
\newtheorem{assumption}{Assumption}
\newtheorem{proposition}{Proposition}
\begin{document}

\begin{frontmatter}

\title{The turnpike property in nonlinear optimal control --- A geometric approach} 

\thanks[footnoteinfo]{This paper was not presented at any IFAC 
meeting. Corresponding author: Noboru Sakamoto. \\
a: This work was partially funded by by JSPS KAKENHI Grant Numbers JP26289128, JP19K04446 and by Nanzan University Pache Research Subsidy I-A-2 for 2019 academic year.\\
b: Supported, in part, 
by the Alexander von Humboldt-Professorship program, 
by the European Research Council (ERC) under the European Union’s Horizon 2020 research and innovation programme (grant agreement No. 694126-DyCon), 
by grant MTM2017-92996 of MINECO (Spain), 
by ELKARTEK project KK-2018/00083 ROAD2DC of the Basque Government,
by ICON of the French ANR and Nonlocal PDEs: Analysis, Control and Beyond, 
and by AFOSR Grant FA9550-18-1-0242.}

\author[NS]{Noboru Sakamoto}\ead{noboru.sakamoto@nanzan-u.ac.jp},    
\author[EZ]{Enrique Zuazua}\ead{enrique.zuazua@fau.de}  

\address[NS]{Faculty of Science and Engineering, Nanzan University, Yamazato-cho 18, Showa-ku, Nagoya, 464-8673, Japan}  
\address[EZ]{Chair in Applied Analysis, Alexander von Humboldt-Professorship, Department of Mathematics, Friedrich-Alexander-Universit\"at Erlangen-N\"urnberg
91058 Erlangen, Germany\\
Departamento de Matem\'aticas, Universidad Aut\'onoma de Madrid, 28049 Madrid, Spain\\
Chair of Computational Mathematics, Fundaci\'on Deusto, Avda Universidades, 24, 48007, Bilbao, Basque Country, Spain}        


\begin{abstract}
This paper presents, using dynamical system theory, a framework for investigating the turnpike property in nonlinear optimal control. First, it is shown that a turnpike-like property appears in general dynamical systems with hyperbolic equilibrium and then, apply it to optimal control problems to obtain sufficient conditions for the turnpike occurs. 
The approach taken is geometric and gives insights for the behaviors of controlled trajectories, allowing us to find simpler proofs for existing results on the turnpike properties. Attempts to remove smallness restrictions  for initial and target states are also discussed based on the geometry of (un)stable manifold and exponential stabilizability of control systems. 
\end{abstract}

\end{frontmatter}
\section{Introduction}
The turnpike property was first recognized in the context of optimal growth by economists (see, e.g., \cite{Mckenzie:63:econometrica}). The turnpike theorems say that for a long-run growth, regardless of starting and ending points, it will pay to get into a growth phase, called {\em von Neumann path}, in the most of intermediate stages. It is exactly like a turnpike and a network of minor roads; "if origin and destination are far enough apart, it will always pay to get on the turnpike and cover distance at the best rate of travel ..."  (quoted from \cite{Dorfman:58:LPEA}). 

In control theory, independently of the turnpike theorems in econometrics, this property was investigated as {\em dichotomy} in linear optimal control \cite{Wilde:72:ieeetac,Rockafellar:73:jota} and later extended to nonlinear systems \cite{Anderson:87:automatica}. In optimal control, the turnpike property essentially means that the solution of an optimal control problem is determined by the system and cost function and independent of time intervals, initial and terminal conditions except in the thin layers at the both ends of the time interval (see, e.g., \cite{Carlson:91:IHOC,Zaslavski:06:TPCVOC,Rapaport:04:esaim}). In the last decades, much progress has been made in the theory of turnpike for finite or infinite dimensional and linear or nonlinear control systems. In \cite{Porretta:13:sicon}, the authors study the turnpike for linear finite and infinite dimensional systems and derive a simple but meaningful inequality, for which we term {\em turnpike inequality}. Their works are extended to finite-dimensional nonlinear systems \cite{Trelat:15:jde}, the semi-linear heat equation \cite{Porretta:16:INdAM}, the wave equation \cite{Gugat:16:syscon,Zuazua:17:arc}, periodic turnpike for systems in Hilbert spaces\cite{Trelat:18:sicon}, optimal shape design \cite{Lance2019_ifacnolcos}, optimal boundary control for hyperbolic systems \cite{Gugat:20:sicon} and general evolution equations \cite{Grune:20:jde}. 
The turnpike property draws attentions of system theory researchers from the viewpoints of model predictive control \cite{Grune:13:automatica,Faulwasser:17:ejc}, dissipative systems (see, e.g., \cite{Willems:71:ac,Willems:72:arma}) \cite{Berberich:18:ieeecst,Damm:14:sicon,Faulwasser:17:automatica,Grune:18:sicon,Grune:16:syscon}, mixed-integer optimal control \cite{Faulwasser:20:ieee_csl}, mechanical systems \cite{Faulwasser:19:nolcos} and the maximum hands-off control \cite{Nagahara:16:ieeetac,sakamoto:20:cdc}.  

In this paper, we first show that turnpike-like behaviors naturally appear in general dynamical systems with hyperbolic equilibrium. The main technique we use is the $\lambda$-lemma which describes trajectory behaviors near invariant manifolds such as stable and unstable manifolds. That the turnpike-like inequality holds implies that if one fixes two ends of a trajectory close to stable and unstable manifolds and designates the time duration sufficiently large, then the trajectory necessarily converges to these manifolds to spend the most of time near the equilibrium. This is exactly the turnpike property. It should be noted that the two ends, as long as they are close to the manifold, do not need to lie in the vicinity of the equilibrium, from which it may be possible to remove locality restrictions in the research of turnpike. 

We apply this inequality to a class of optimal control problems in which terminal states are not specified and the steady optimal solutions are not the origin as in \cite{Porretta:13:sicon,Porretta:16:INdAM,Zuazua:17:arc,Trelat:18:sicon} as well as to a class of optimal control problems in which two terminal states are specified and the steady optimal point is the origin as in \cite{Wilde:72:ieeetac,Anderson:87:automatica,Trelat:15:jde}. For both classes of problems, we employ a Dynamic Programming approach with Hamilton-Jacobi equations (HJEs). The characteristic equations for HJEs are Hamiltonian systems and the stabilizability (controllability for the second class) and detectability conditions assure that the equilibrium of the Hamiltonian systems is hyperbolic. The {\em controlled trajectories} appear in these Hamiltonian systems and we apply the turnpike result for dynamical system. Then, the existence of the trajectory satisfying initial and boundary conditions is guaranteed. In this paper, we derive {\em sufficient conditions} for optimality by using the Dynamic Programming and HJEs and by imposing a condition that guarantees the existence of the solution to the HJEs (Lagrangian submanifold property, see, e.g., \cite{Libermann:87:SGAM}). 

The present manuscript expands upon our conference contribution \cite{sakamoto:19:cdc} incorporating a new result on the relationship between nonlinear stabilizability and existence of infinite horizon optimal control \cite{Sakamoto:20:prep}. It allows one to give an estimate of the existence region of a stable manifold of hyperbolic Hamiltonian system associated with an optimal control problem, from which one may be able to predict the occurrence of turnpike (see Section~\ref{sctn:ex_2nd_class}). The manuscript also contains a number of examples worked out to show how the proposed geometric approach is effectively applied for turnpike analysis and appendices for necessary results in the theory of algebraic Riccati equations and for stable manifold estimate in Hamiltonian systems. 

The organization of the paper is as follows. In Section~\ref{sctn:lamda_lemma} we review some of key tools from dynamical system theory and derive the turnpike inequality. In Section~\ref{sctn:ocp}, we apply it to optimal control problems. Section~\ref{sbsctn:ocp1} handles the problem where boundary state is free and Section~\ref{sbsctn:opc2} handles the problem where initial and boundary states are fixed. Section~\ref{sctn:examples} shows turnpike analyses for a class of nonlinear systems for which target $z$ in $(\mathrm{OCP}_1)$ can be taken arbitrarily large and a class of nonlinear systems for which initial states can be taken arbitrarily large. Section~\ref{sctn:disscussion} discusses possible extensions for more general turnpike using the geometric approach.  

\section{Turnpike in dynamical systems}\label{sctn:lamda_lemma}
Let us consider a nonlinear dynamical system of the form
\begin{equation}
    \dot z= f(z),\label{eqn:dyn_sys}
\end{equation}
where $f:\mathbb{R}^N\to\mathbb{R}^N$ is of $C^r$ class ($r\geqslant1$). We assume that $f(0)=0$ and {\em the hyperbolicity of $f$ at $0$}, namely, assume that  $Df(0)\in\mathbb{R}^{N\times N}$ has $k$ eigenvalues with strictly negative real parts and $N-k$ eigenvalues with strictly positive real parts. 

It is known, as {\em the stable manifold theorem}
, that there exist $C^r$ manifolds $S$ and $U$, called {\em stable manifold} and {\em unstable manifold} of (\ref{eqn:dyn_sys}) at $0$, respectively, defined by 
\begin{align*}
S&:=\{z\in\mathbb{R}^N\,|\, \varphi(t,z)\to0 \text{ as } t\to \infty\},\\
U&:=\{z\in\mathbb{R}^N\,|\, \varphi(t,z)\to0 \text{ as } t\to -\infty\},
\end{align*}
where $\varphi(t,z)$ is the solution of (\ref{eqn:dyn_sys}) starting $z$ at $t=0$. Let $E^s$, $E^u$ be stable and unstable subspaces in $\mathbb{R}^N$ of $Df(0)$ with dimension $k$, $N-k$, respectively. It is known that $S$, $U$ are invariant under the flow of $f$ and are tangent to $E^s$, $E^u$, respectively, at $z=0$. See, e.g., \cite{Hale:73:ODE,Palis:82:GTDS} for more details on the theory of stable manifold. 

We will consider limiting behavior of submanifolds under the flow of $f$ and need to introduce topology for maps and manifolds. Let $M$ be a compact manifold of dimension $m$ and the space $C^r(M,\mathbb{R}^l)$ of $C^r$ maps, $0\leqslant r<\infty$, defined on $M$. There exists a natural vector space structure on $C^r(M,\mathbb{R}^l)$. Since $M$ is compact, we take a finite cover of $M$ by open sets $V_1,\ldots,V_k$ and take a local chart $(z_i,U_i)$ for $M$ with $z_i(U_i)=B(2)$ such that $z_i(V_i)=B(1)$, $i=1,\ldots,k$, where $B(1)$ and $B(2)$ are the balls of radius 1 and 2 at the origin of $\mathbb{R}^m$. For a map $g\in C^r(M,\mathbb{R}^l)$, we define a norm by
\begin{multline*}
\|g\|_r:=\max_i \sup\{|g(u)|,\|Dg^i(u)\|, \ldots, \\
    \|D^rg^i(u)\|\,|\,u\in B(1)\},
\end{multline*}
where $g^i=g\circ {z_i}^{-1}$, local representation of $g$, and $\|\cdot\|$ is a norm for linear maps. It is known that $\|\cdot\|_r$ does not depend on the choice of finite cover and we call it $C^r$ norm. For maps in $C^r(M,N)$ where $N$ is a manifold, we embed $N$ in a Euclidean space with sufficiently high dimension. Let $L$, $L'$ be $C^r$ submanifolds of $M$ and let $\varepsilon>0$. We say that {\em $L$ and $L'$ are $\varepsilon$ $C^r$-close} if there exists a $C^r$ diffeomorphism $h:L\to L'$ such that $\|i'\circ h-i\|_r<\varepsilon$, where $i:L\to M$ and $i':L'\to M$ are inclusion maps. In this case, we use the notation $d_L^r(L'):=\|i'\circ h-i\|_r$. 

By a $k$-dimensional (topological) disc we mean a set that is homeomorphic to $D^k:=\{ (x_1,\ldots,x_k)\in\mathbb{R}^k\,|\, x_1^2+\cdots+x_k^2\leqslant1 \}$. The following lemma is known as {\em the $\lambda$-lemma} or {\em inclination lemma} and plays a crucial role in the theory of dynamical systems (see \cite{Palis:82:GTDS,Wen:16:DDS}). 
\begin{lem}[The $\lambda$-lemma] Suppose that $z=0$ is a hyperbolic equilibrium for (\ref{eqn:dyn_sys}). Suppose also that $S$ and $U$ are $k$, $(N-k)$-dimensional
stable and unstable manifolds of $f$ at $0$, respectively. For any $(N-k)$-dimensional disc $B$ in $U$, any point $z\in S$, any $(N-k)$-dimensional disc $D$ transversal to $S$ at $z$ and any $\varepsilon>0$, there exists a $T>0$ such that if $t>T$, $\varphi(t,D)$ contains an $(N-k)$-dimensional disc $\tilde D$ with $d_{B}^1(\tilde D)<\varepsilon$. 
\end{lem}

\begin{figure}[htb]
    \centering
    \includegraphics[width=0.3\textwidth]{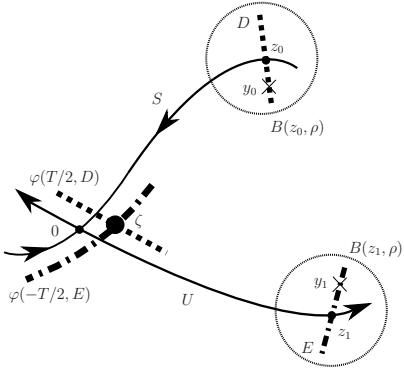}
    \caption{A scheme of the proof of Theorem~\ref{thm:turnpk_in_dyn}}
    \label{fig:turnpk_in_dyn}
\end{figure}

Next, we show that {\em the turnpike behavior} appears in the transition of points near the stable manifold to points near the unstable manifold if the transition duration is designated large. Let $z_0\in S$ and $z_1\in U$ be given points. From the stable manifold theorem, it holds that 
\begin{subequations}\label{ineq:s_mani_u_mani}
\begin{align}
    |\varphi(t,z_0)| & \leqslant K e^{-\mu t} \text{ for } t\geqslant0,\\
    |\varphi(t,z_1)| & \leqslant K e^{\mu t} \text{ for } t\leqslant0, 
\end{align}
\end{subequations}
where $K>0$ is a constant dependent on $z_0$ and $z_1$ and $\mu>0$ is a constant independent of $z_0$ and $z_1$. 
\begin{prop}\label{prop:pre-turnpike}
\begin{enumerate}[(i)]
    \item There exists a $T_0>0$ such that for every $T>T_0$ there exists a $\rho=\rho(T)>0$ such that 
    \[
    |\varphi(t,y)|\leqslant Ke^{-\mu t} \text{ for } t\in[0,T],\ y\in B(z_0,\rho), 
    \]
    where $B(z_0,\rho)$ is the $N$-dimensional open ball centered at $z_0$ with radius $\rho$. Moreover, $\rho\to0$ when $T\to\infty$.
\item There exist a $T_0<0$ such that for every $T<T_0$ there exists a $\rho=\rho(T)>0$ such that 
    \[
    |\varphi(t,y)|\leqslant Ke^{\mu t} \text{ for } t\in[T,0],\ y\in B(z_1,\rho).
    \]
Moreover, $\rho\to0$ when $T\to -\infty$.
\item For any $(N-k)$-dimensional disc $\bar D$
transversal to $S$ at $z_0$ and any $k$-dimensional disc $\bar E$ transversal to $U$ at $z_1$, there exists a $T_0>0$ such that for any $T>T_0$ there exist an $(N-k)$-dimensional disc $D\subset \bar D$ transversal to $S$ at $z_0$ and a $k$-dimensional disc $E\subset \bar E$ transversal to $U$ at $z_1$ such that $\varphi(T,D)$ intersects $\varphi(-T,E)$ at a single point. 
\end{enumerate}
\end{prop}
%
\quad{\em Proof.} 
(i), (ii) These are consequences of the properties (\ref{ineq:s_mani_u_mani}) of (un)stable manifold, which can be derived using contradiction arguments.  The proofs of the convergence $\rho\to0$ as $T\to\pm\infty$ also use contradictions to the facts that $S$, $U$ are submanifolds with strictly lower dimension than $N$.\\
(iii) First we take an $(N-k)$-dimensional disc $U_0$ in $U$ passing through 0, a $k$-dimensional disc $S_0$ in $S$ passing through 0 and an $\varepsilon>0$ arbitrarily. From the $\lambda$-lemma, there exists a $T_0>0$ such that for any $T>T_0$ there exists an $(N-k)$-dimensional disc $D\subset \bar D$ transversal to $S$ at $z_0$ and a $k$-dimensional disc $E\subset \bar E$ transversal to $U$ at $z_1$ such that $d_{U_0}^1(\varphi(T,D))<\varepsilon$, $d_{S_0}^1(\varphi(-T,E))<\varepsilon$. Since $E^s\cap E^u=\{0\}$, it is possible to take $\varepsilon$, $S_0$ and $U_0$ so that $\varphi(T,D)$ intersects $\varphi(-T,E)$ at a single point. $\hfill\blacksquare$
%
\begin{rem}
It should be noted that the above statements, especially (i) and (ii), are only on finite interval $[0,T]$. This is the major difference from the trajectories on the stable and unstable manifolds.
\end{rem}
\begin{thm}\label{thm:turnpk_in_dyn}
For any $z_0\in S$, any $z_1\in U$, any $(N-k)$-dimensional disc $\bar D$
transversal to $S$ at $z_0$ and any $k$-dimensional disc $\bar E$ transversal to $U$ at $z_1$, there exists a $T_0>0$ such that for every $T>T_0$ there exist $\rho=\rho(T)>0$, $y_0\in B(z_0,\rho)\cap \bar D$ and $y_1\in B(z_1,\rho)\cap \bar E$ such that $\varphi(T,y_0)=y_1$ and 
\[
|\varphi(t,y_0)|\leqslant K\left[e^{-\mu t}+e^{-\mu(T-t)}\right] \text{ for } t\in [0,T].
\]
Moreover, $\rho\to0$ when $T\to\infty$. 
\end{thm}
\quad{\em Proof.} 
Take the largest $T_0$ and the smallest $\rho$ in Proposition~\ref{prop:pre-turnpike}. We rename this $T_0$ as $T_0/2$. Take arbitrary $T>T_0$ and use Proposition~\ref{prop:pre-turnpike}-(iii) to get a disc $D$ which is $(N-k)$-dimensional and transversal to $S$ at $z_0$ and a disc $E$ which is $k$-dimensional and transversal to $U$ at $z_1$ satisfying $D\subset B(z_0,\rho)$ and $E\subset B(z_1,\rho)$. This is possible by taking smaller $S_0$ and $U_0$ in the proof of Proposition~\ref{prop:pre-turnpike}-(iii). Then, there exists a single point $\zeta$ such that $\varphi(T/2,D)\cap\varphi(-T/2,E)=\{\zeta\}$ (see Fig~\ref{fig:turnpk_in_dyn}). Let $y_0:=\varphi(-T/2,\zeta)$. Then, $y_0\in D\subset B(z_0,\rho)$ and by Proposition~\ref{prop:pre-turnpike}-(i), we have
\begin{equation}\label{eqn:est_stable}
    |\varphi(t,y_0)|\leqslant Ke^{-\mu t} \text{ for }0\leqslant t\leqslant T/2. 
\end{equation}
Let $y_1:=\varphi(T/2,z)$. Then, $y_1\in E\subset B(x_1,\rho)$ and 
\[
    |\varphi(t,y_1)|\leqslant Ke^{\mu t}\text{ for }-T/2\leqslant t\leqslant0.
\]
This shows that 
\[
    |\varphi(t+T,y_0)|\leqslant Ke^{\mu t}\text{ for }-T/2\leqslant t\leqslant0, 
\]
and consequently, 
\begin{equation}\label{eqn:est_unstable}
    |\varphi(t,y_0)|\leqslant Ke^{-\mu(T-t)}\text{ for } T/2\leqslant t\leqslant T.
\end{equation}
Combining (\ref{eqn:est_stable}) and (\ref{eqn:est_unstable}), we get the inequality in the theorem. The last assertion follows from Proposition~\ref{prop:pre-turnpike}-(i) and (ii). $\hfill\blacksquare$
\section{Turnpike in nonlinear optimal control}\label{sctn:ocp}
Let us consider a nonlinear control system
\begin{equation}
\dot x= f(x)+g(x)u,\ x(t_0)=x_0,\label{eqn:nsys_general}
\end{equation}
where $f:\mathbb{R}^n\to\mathbb{R}^n$, $g:\mathbb{R}^n\to\mathbb{R}^{n\times m}$ are of $C^2$ class with $f(0)=0$, $x(t)\in\mathbb{R}^n$ is state variables and $u(t)\in\mathbb{R}^m$ is control input. An optimal control problem or OCP is to find a control input for (\ref{eqn:nsys_general}) such that the cost functional
\[
J(u)=\int_{0}^TL(x(t),u(t))\,dt
\]
is minimized, where we set $J(u)=+\infty$ when the existence domain of solution for (\ref{eqn:nsys_general}) is strictly contained in $[0,T)$. There are several types in OCPs depending on whether or not the terminal time $T$ is specified and whether or not the state variables are specified at the terminal time. In this paper, we consider OCPs where the terminal time $T$ is specified and two types of OCPs; one in which the state variables are free at $t=T$ and another in which they are fixed at $t=T$. 
For both types of OCPs, we are interested in the relationship between the solution $u_T$ and corresponding trajectory $x_T$ of an OCP and steady-state optimum pair $(\bar u, \bar x)$, which will be defined more precisely later on. 
\begin{defn}\cite{Carlson:91:IHOC}\label{dfn:turnpike}
An optimal pair $(u_T,x_T)$ has the turnpike property if for any $\varepsilon>0$, there exists an $\eta_\varepsilon>0$ such that 
\[
\left| \left\{t\geqslant0\,|\,|u_T(t)-\bar{u}|+|x_T(t,x_0)-\bar{x}|>\varepsilon  \right\} \right| <\eta_\varepsilon
\]
for all $T>0$, where $\eta_\varepsilon$ depends only on $\varepsilon$, $f$, $g$, $x_0$, and $L$ and $|\cdot|$ denotes length (Lebesgue measure) of interval. 
\end{defn}
\begin{rem}{\rm
{\em Turnpike inequality} is to require $x_T$ and $u_T$ to satisfy 
\begin{equation}
    |u_T (t)-\bar{u}|+|x_T(t,x_0)-\bar{x}|\leqslant 
    K\left[e^{-\mu t}+e^{-\mu(T-t)}\right]\label{eqn:turnpike_ineq_control}
\end{equation}
for some constants $K>0$ and $\mu>0$ independent of $T$, which is a sufficient condition for the turnpike property in Definition~\ref{dfn:turnpike}. Also, it should be noted that requiring (\ref{eqn:turnpike_ineq_control}) limits ourselves to {\em the exponential input-state turnpike} defined in \cite{Grune:16:syscon}. 
}
\end{rem}

\subsection{The OCP with state variables unspecified at the terminal time}\label{sbsctn:ocp1}
For system (\ref{eqn:nsys_general}), we consider the following cost functional
\[
J_1(u)=\frac{1}{2}\int_{0}^T|Cx(t)-z|^2+|u(t)|^2\, dt,
\]
where $C\in\mathbb{R}^{r\times n}$ and $z\in\mathbb{R}^r$ is a given vector ({\em target)}. We call this problem $(\text{OCP}_1)_T$; 
\begin{align*}
(\text{OCP}_1)_T:\ \  &\text{Find a control }u\in L^\infty(0,T;\mathbb{R}^m) \text{ such that }\\
& J_1(u) \text{ along } (\ref{eqn:nsys_general}) \text{ is minimized}\text{ over all } \\
&u\in L^\infty(0,T,\mathbb{R}^m).
\end{align*}
Associated with $(\text{OCP}_1)_T$, we consider a steady optimization problem
\begin{align*}
\text{(SOP)}: \quad &\text{Minimize }J_s(x,u) = \frac{1}{2}(|Cx-z|^2+|u|^2) \\
&\text{ over all }(x,u)\in\mathbb{R}^n\times\mathbb{R}^m \text{such that }\\
&f(x)+g(x)u=0. 
\end{align*}
We assume the following. 
\begin{assum}\label{assm:steady}
(SOP) has a solution $(\bar{x},\bar{u})=(\bar{x}(z),\bar{u}(z))$. 
\end{assum}
\vskip 1ex
Also, associated with $(\text{OCP}_1)_T$, we can derive a Hamilton-Jacobi equation
\begin{gather}
\begin{multlined}
    V_t(t,x)+V_x(t,x)f(x)\\
    -\frac{1}{2}V_x(t,x)g(x)g(x)^\top V_x(t,x)^\top +\frac{1}{2}|Cx-z|^2=0,
\end{multlined}\label{eqn:hje}\\
    V_x(T,x)=0,\label{eqn:hjc}
\end{gather}
for $V(t,x)$, where $V_t=D_tV$, $V_x=D_xV$. Defining a Hamiltonian 
\[
H(x,p)=p^\top f(x)-\frac{1}{2}p^\top g(x)g(x)^\top p+\frac{1}{2}|Cx-z|^2, 
\]
we consider the corresponding characteristic equation for (\ref{eqn:hje})-(\ref{eqn:hjc})
\begin{equation}
\dot x_i=\frac{\partial H}{\partial p_i},\ \dot p_i=-\frac{\partial H}{\partial x_i},\quad i=1,\ldots,n \label{eqn:Hsys}
\end{equation}
with $\ p_i(T)=0$, $i=1,\ldots,n.$ 
Note that since the system (\ref{eqn:nsys_general}) is time-invariant, the equation corresponding to $V_t$ is not necessary. The following fact is readily verified.\\[1ex]
{\bf Fact.} A solution $(\bar{x},\bar{u})$ of (SOP) corresponds to an equilibrium point $(\bar{x},\bar{p})$ of (\ref{eqn:Hsys}) with $\bar{u}=-g(\bar{x})^\top \bar{p}$.

Let $A_z=D_xD_pH(\bar{x},\bar{p})$, $B_z=g(\bar{x})$. 
\begin{assum}\label{assm:stb_detc}
The triplet $(C,A_z,B_z)$ is stabilizable and detectable. 
\end{assum}
\vskip 1ex
Under Assumptions~\ref{assm:steady}, \ref{assm:stb_detc}, the equilibrium $(\bar{x},\bar{p})$ is hyperbolic equilibrium for the Hamiltonian system (\ref{eqn:Hsys}) and there exist stable and unstable manifolds for (\ref{eqn:Hsys}) at $(\bar{x},\bar{p})$ which are expressed as 
\begin{equation}
S_z=\Tilde{S}+\{(\bar{x},\bar{p})\},\quad U_z=\Tilde{U}+\{(\bar{x},\bar{p})\}.\label{eqn:s-u_manifolds}
\end{equation}
Here, $\tilde{S}$, $\tilde{U}$ are the stable and unstable manifold of (\ref{eqn:Hsys}) in the coordinates $(\tilde{x},\tilde{p})$, where $\tilde{x}=x-\bar{x}$, $\tilde{p}=p-\bar{p}$, which is re-written as  
\[
\frac{d}{dt}\bmat{\tilde{x}\\\tilde{p}}=\bmat{A_z&-B_zB_z^\top \\-C^\top C&-A_z^\top }\bmat{\tilde{x}\\\tilde{p}}+o(|\tilde{x}|+|\tilde{p}|).
\]

We can now state the main theorem of this subsection. Let $\pi_1:(x,p)\mapsto x$, $\pi_2:(x,p)\mapsto p$ be canonical projections. 
\begin{thm}\label{thm:trunpike_control_type1}
Under Assumptions~\ref{assm:steady}, \ref{assm:stb_detc}, suppose that $x_0\in \mathrm{Int}(\pi_1(S_z))$, where $\mathrm{Int}(\cdot)$ is the interior of a set in $\mathbb{R}^n$, and that $U_z$ intersects $p=0$ transversally. If $T$ is taken sufficiently large, then there exists a solution $(x_T(t,x_0),p_T(t,x_0))$ to (\ref{eqn:Hsys}) satisfying $x_T(0,x_0)=x_0$ and $p_T(T,x_0)=0$. If, moreover, 
\begin{equation}
\det D_{x_0}x_T(t,x_0)\ne 0 \text{ for } t\in[0,T],\label{eqn:p-condition}
\end{equation}
then 
\[
u_T (t):=-g(x_T(t,x_0))^\top p_T(t,x_0)
\]
is the local optimal solution for $(\mathrm{OCP}_1)_T$ and turnpike inequality (\ref{eqn:turnpike_ineq_control}) holds for some constants $K>0$ and $\mu>0$ which are independent of $T$. 
\end{thm}
\quad{\em Proof.} 
Let $\zeta_0=(x_0,0)$ and $\zeta_1=(x_1,0)$, where $(x_1,0)\in U_z$ and take $n$-dimensional discs $\{(x_0,p)\,|\,|p|<\rho\}$, $\{(x,0)\,|\,|x-x_1|<\rho\}$, 
which correspond to $z_0$, $z_1$, $B(z_0,\rho)\cap\bar D$ and $B(z_1,\rho)\cap\bar E$ in Theorem~\ref{thm:turnpk_in_dyn}, respectively. Then the theorem implies that for a sufficiently large $T>0$, there exist $\rho=\rho(T)>0$, $p_0$ and $x_1'$ with $|p_0|<\rho$, $|x_1'-x_1|<\rho$ such that 
\[
\varphi(T,(x_0,p_0))=(x'_1,0), 
\]
where $\varphi(t,(x_0,p_0))$ denotes the solution of (\ref{eqn:Hsys}) starting from $(x_0,p_0)$. This shows that the two-point boundary value problem associated with $(\mathrm{OCP}_1)_T$ has been solved. Let $(x_T(t,x_0), p_T(t,x_0))=\varphi(t, (x_0,p_0))$. Then, the theorem also says that  
there exist $K'>0$ and $\mu>0$ such that 
\begin{multline*}
|x_T(t,x_0)-\bar{x}|+|p_T(t,x_0)-\bar{p}|\leqslant K'[e^{-\mu t}+e^{-\mu(T-t)}]\\
\text{ for }0\leqslant t\leqslant T.
\end{multline*}
Since $|u_T(t)-\bar{u}| \leqslant \sup \|g(x)\||p(t,x_0)-\bar{p}|$, (\ref{eqn:turnpike_ineq_control}) holds with $K=2K'(1+\sup \|g(x)\|)$, where supremum is taken along the trajectory. The condition (\ref{eqn:p-condition}) guarantees that there exists a Lagrangian submanifold in a neighborhood of this trajectory and this implies the existence of solution $V(t,x)$ to (\ref{eqn:hje})-(\ref{eqn:hjc}) in the neighborhood. Then, the verification theorem in Dynamic Programming (see, e.g., \cite{Athans:66:OC}) shows that the control $u^\ast$ is locally optimal. $\hfill\blacksquare$
\begin{rem}{\rm
The condition (\ref{eqn:p-condition}) guarantees that the solution $V$ to (\ref{eqn:hje}) exists in a neighborhood of the trajectory $(x_T(t,x_0),p_T(t,x_0))$, $0\leqslant t\leqslant T$. The optimality of $u_T$ is valid only in the neighborhood. This existence theory is described using the notion of {\em Lagrangian submanifold} (see, e.g., \cite{Libermann:87:SGAM}) and when one seeks for larger domain of existence, the non-uniqueness issue of solution arises. We refer to \cite{Day:98:math-syst-estim-contr} for general analysis of non-unique solutions and \cite{Osinga:06:jdde,Horibe:17:ieee_cst,Horibe:19:ieee_cst} for non-unique optimal controls for mechanical systems.  
}
\end{rem}
We next show that for small $x_0$, $z$, the solution for $(\mathrm{OCP}_1)_T$ has a solution with turnpike property using perturbation theory of stable manifold. Let $A=f(0)$, $B=g(0)$. 
\begin{assum}\label{assm:stabilizable_at_0}
The triplet $(C,A,B)$ is stabilzable and detectable.
\end{assum}
{\bf Fact.} Under Assumption~\ref{assm:stabilizable_at_0}, there is a neighborhood of $z=0$ in $\mathbb{R}^r$ such that (SOP) has a unique solution for $z$ in the neighborhood and $(C,A_z,B_z)$ is stabilzable and detectable.
\begin{cor}\label{cor:turnpike_nonlinear_control_1}
Under Assumption~\ref{assm:stabilizable_at_0}, for sufficiently small $x_0$ and $z$ and for sufficiently large $T$, $(\mathrm{OCP}_1)_T$ has a solution with the turnpike property. 
\end{cor}
\quad {\em Proof.} From the Fact above, under Assumption~\ref{assm:stabilizable_at_0}, the Hamiltonian system (\ref{eqn:Hsys}) has stable manifold $S_z$ and unstable manifold $U_z$ at $(\bar{x},\bar{p})$. For $z=0$ the linear part of the Hamiltonian system is 
$\mathrm{Ham}=\left[\begin{smallmatrix}A&-BB^\top \\-C^\top C&-A^\top \end{smallmatrix}\right]$, for which we apply the eigen structure analysis in Appendix~\ref{sctn:appx}.  
Apply Lemma~\ref{lemma:ham_appdx} with $R=BB^\top$, $Q=C^\top C$ and let $P$ and $L$ as in the Appendix. Then, the tangent spaces $T_0S_0$, $T_0U_0$ of $S_0$, $U_0$ at the origin can be written as
\begin{align*}
    T_0S_0&=\{(u,Pu)\,|\,u\in\mathbb{R}^n \},\\
    T_0U_0&=\{ (Lu,(PL+I)u)\,|\, u\in\mathbb{R}^n \}. 
\end{align*}
 From the expression of $T_0S_0$, one can take $x_0$ sufficiently small so that there is an $n$-dimensional disc $D_0$ in $S_0$ that contains the origin and $x_0$ in its interior. 
From Lemma~\ref{lemma:PL+I}, 
$PL+I$ is nonsingular and therefore, $T_0U_0$ intersects $p=0$ transversally, which implies that there is an $n$-dimensional disc $E_0$ in $U_0$ that intersects $p=0$ transversally. Let $X_H(x,p;z)$ be the Hamiltonian vector field (\ref{eqn:Hsys}). As $z\to0$, $X_H(x,p,z)$ can be arbitrarily close to $X_H(x,p;0)$ with $C^1$ topology in an appropriate compact set. The stable manifold theory (see, e.g., \cite[Theorem 6.2]{Palis:82:GTDS}) ensures that there exists a small $z$ so that there are $n$-dimensional discs $D_z\subset S_z$, $E_z\subset U_z$ that are close enough to $D_0$, $E_0$, respectively, with $C^1$-topology. For this $z$, it holds that $x_0\in\mathrm{Int}(\pi_1(D_z))$ and $E_z$ intersects $p=0$ transversally. Now, all the hypotheses in Theorem~\ref{thm:trunpike_control_type1} are satisfied. $\hfill\blacksquare$

Next Corollary is proved in \cite{Porretta:13:sicon,Trelat:18:sicon} in the study of the turnpike property for infinite dimensional systems under slightly more restrictive conditions (controllability and observability rather than stabilizability and detectability). Their proofs are based on the estimates on adjoint variables in the linear Hamiltonian system (\ref{eqn:Hsys}) which is derived as a necessary condition of optimality. Here, we give an alternative proof using the geometric picture in Theorem~\ref{thm:trunpike_control_type1}. 
\begin{cor}\label{cor:linear_turnpike}
Suppose that the system (\ref{eqn:nsys_general}) is linear, that is, $f(x)=Ax$ and $g(x)=B$ with real constant matrices $A\in\mathbb{R}^{n\times n}$ and $B\in\mathbb{R}^{n\times m}$. Under Assumption~\ref{assm:stabilizable_at_0}, $(\mathrm{OCP}_1)_T$ has the global solution $u^\ast(t)$, $0\leqslant t\leqslant T$ for any $z\in\mathbb{R}^r$. Moreover, turnpike inequality (\ref{eqn:turnpike_ineq_control}) holds.  
\end{cor}
\quad {\em Proof.} We use some of the notations in the proof of  Corollary~\ref{cor:turnpike_nonlinear_control_1}. 
The unique solution $(\bar{x},\bar{p})$ to (SOC) is expressed as $\left[\begin{smallmatrix}\bar{x}\\\bar{p} \end{smallmatrix}\right]=-\mathrm{Ham}^{-1}\left[\begin{smallmatrix}0\\C^\top z \end{smallmatrix}\right]$. $U_z$ and $S_z$ in (\ref{eqn:s-u_manifolds}) can be written as
\begin{align*}
    S_z&=\{(u,Pu)\,|\,u\in\mathbb{R}^n \}+\{(\bar{x},\bar{p})\},\\
    U_z&=\{ (Lu,(PL+I)u)\,|\, u\in\mathbb{R}^n \}+\{(\bar{x},\bar{p})\}. 
\end{align*}
It is readily seen that $x_0\in \mathrm{Int}(\pi_1(S))$ for any $x_0\in\mathbb{R}^n$ and
$U$ intersects $p=0$ transversally for any $z\in\mathbb{R}^r$. The condition (\ref{eqn:p-condition}) is equivalent to the nonsingularity of (1,1)-block in $\exp[t \mathrm{Ham}]$, which is proved in Lemma~\ref{lemma:Phi_11_appdx}. $\hfill\blacksquare$
\begin{rem}{\rm
\begin{enumerate}
\item
Although the problem in Corollary~\ref{cor:linear_turnpike} is linear, it is not an easy task to explicitly write down the solution for (\ref{eqn:hje})-(\ref{eqn:hjc}) except for $z=0$. This corollary, however, says that the solution globally exists. 
\item
As is discussed in \cite{Porretta:16:INdAM,Trelat:15:jde,Phigin:20:semilnear_arxiv}, relaxing the smallness conditions in Corollary~\ref{cor:turnpike_nonlinear_control_1} is one of major challenges in the research of nonlinear turnpike. In \S~\ref{sctn:ex_1st_class}, we show a class of nonlinear OCPs for which turnpike occurs for all $z$ by explicitly analyzing unstable manifold. 
\end{enumerate}
}
\end{rem}

\subsection{The OCP with state variables specified at the terminal time}\label{sbsctn:opc2}
In this subsection, we consider an OCP for (\ref{eqn:nsys_general}) with arbitrarily specified terminal states. Let $x_f\in\mathbb{R}^n$ be given. Let us define cost functional
\[
J_2(u)=\frac{1}{2}\int_0^{T} x(t)^\top C^\top Cx(t)+|u|^2\,dt,
\]
and consider 
\begin{align*}
(\mathrm{OCP}_2)_T:\ \  &\text{Find a control }u\in L^\infty(0,T;\mathbb{R}^m) \text{ such that }\\
& J_2(u) \text{ along (\ref{eqn:nsys_general}) is minimized over all } \\
& u\in L^\infty (0,T;\mathbb{R}^m) \text{ such that }x(T)=x_f.
\end{align*}
With Assumption~\ref{assm:stabilizable_at_0}, the corresponding steady optimization problem has a unique solution $(\bar{x},\bar{u})=(0,0)$ around the origin. 
The Hamilton-Jacobi equation associated with $(\mathrm{OCP}_2)_T$ is 
\begin{gather}
\begin{multlined}
    V_t(t,x)+V_x(t,x)f(x)\\
    -\frac{1}{2}V_x(t,x)g(x)g(x)^\top V_x(t,x)^\top +\frac{1}{2}x^\top C^\top Cx=0. 
\end{multlined}\label{eqn:hje_2}
\end{gather}
The Hamiltonian in this case is
\[
H(x,p)=p^\top f(x)-\frac{1}{2}p^\top g(x)g(x)^\top p+\frac{1}{2}x^\top C^\top Cx, 
\]
and the corresponding characteristic equation for (\ref{eqn:hje_2}) is 
\begin{equation}
\dot x_i=\frac{\partial H}{\partial p_i},\ \dot p_i=-\frac{\partial H}{\partial x_i},\quad i=1,\ldots,n \label{eqn:Hsys_2}
\end{equation}
with $x(0)=x_0$ and $x(T)=x_f$. 

Under Assumptions~\ref{assm:stabilizable_at_0}, the Hamiltonian system (\ref{eqn:Hsys_2}) can be written as
\[
\frac{d}{dt}\bmat{x\\p}=\mathrm{Ham}\bmat{x\\p}+\text{higher oder terms}, 
\]
and 
the origin is a hyperbolic equilibrium with $n$ stable and $n$ unstable eigenvalues. Let $S$ and $U$ be the stable and unstable manifolds of (\ref{eqn:Hsys_2}) at the origin. 
\begin{thm}\label{thm:trunpike_control_type2}
Under Assumption~\ref{assm:stabilizable_at_0}, suppose that $x_0\in\mathrm{Int}(\pi_1(S))$ and $x_f\in\mathrm{Int}(\pi_1(U))$. If $T>0$ is taken sufficiently large, there exists a solution $(x_T(t,x_0),p_T(t,x_0))$ to (\ref{eqn:Hsys_2}) satisfying $x(0)=x_0$ and $x(T)=x_1$. If, moreover, 
\begin{equation}
\det D_{x_0}x(t,x_0)\ne 0 \text{ for } t\in[0,T],\label{eqn:p-condition_2}
\end{equation}
then
\[
u_T (t)=-g(x_T(t,x_0))^\top p_T(t,x_0)
\]
is the local optimal solution for $(\mathrm{OCP}_2)_T$ and turnpike inequality (\ref{eqn:turnpike_ineq_control}) hols for some $K>0$, $\mu>0$ independent of $T$. 
\end{thm}
\quad{\em Proof.} Let $\zeta_0=(x_0,0)$, $\zeta_1=(x_f,0)$, $\{(x_0,p)\,|\,|p|<\rho\}$ and $\{(x_f,p)\,|\,|p|<\rho\}$ which correspond to $z_0$, $z_1$, $B(z_0,\rho)\cap \bar{D}$ and $B(z_1,\rho)\cap\bar{E}$ 
in Theorem~\ref{thm:turnpk_in_dyn}, respectively. 
Then, for $T>0$ sufficiently large, there exist $\rho>0$, $p_0$ and $p_1$ with $|p_0|,|p_1|<\rho$ 
 such that a solution to (\ref{eqn:Hsys_2}) connecting $(x_0,p_0)$ and $(x_f,p_1)$ exists. The rest of the proof is almost the same as Theorem~\ref{thm:trunpike_control_type1}.$\hfill\blacksquare$
\begin{cor}\label{cor:tp_under_smallness}
Let us additionally impose the controllability of $(A,B)$ in Assumption~\ref{assm:stabilizable_at_0}. Then, for sufficiently small $|x_0|$ and $|x_f|$ and sufficiently large $T$, the local optimal control exists and turnpike inequality (\ref{eqn:turnpike_ineq_control}) holds. 
\end{cor}
\quad{\em Proof.} 
We again employ the eigen structure analysis (\ref{eqn:eigen_eq_appdx}). The tangent spaces of $S$ and $U$ at the origin are written as 
\begin{align*}
T_0S&=\{ (u,Pu)\,|\, u\in\mathbb{R}^n \},\\
T_0U&=\{(u,(PL+I)L^{-1}u)\,|\, u\in\mathbb{R}^n\}. 
\end{align*}
The latter is obtained by showing, using the controllability of $(A,B)$, that $L$ is strictly negative definite (Lemma~\ref{lemma:PL+I}). Therefore, $x_0\in \mathrm{Int}(\pi_1(S))$ and $x_f\in\mathrm{Int}(\pi_1(U))$ for sufficiently small $|x_0|$, $|x_f|$. It is seen that the condition (\ref{eqn:p-condition_2}) holds for these  $|x_0|$, $|x_f|$ (making them smaller if necessary) from the analysis on $\Phi_{11}(t)$ in the proof of Theorem~\ref{thm:trunpike_control_type1}.$\hfill\blacksquare$
\begin{rem}{\rm
\begin{enumerate}[(i)]
    \item 
The linear counterpart of Corollary~\ref{cor:tp_under_smallness} is in \cite{Wilde:72:ieeetac} where they use anti-stabilizing solution $P_u$ for the Riccati equation. In this case, the turnpike holds for all $x_0$ and $x_f$. It can be shown that $P_u=(PL+I)L^{-1}$. Note that in Corollary~\ref{cor:tp_under_smallness} we only need the detectability condition. Corollary~\ref{cor:tp_under_smallness} is obtained in \cite{Anderson:87:automatica} using Hamilton-Jacobi theory under unusual nonlinear controllability and observability conditions. Compared with their conditions, we use only the {\em linear} controllability and detectability which can be easily checked. The authors of \cite{Trelat:15:jde} obtain similar results to Corollary~\ref{cor:tp_under_smallness} with more generalized terminal conditions. 
\item
Corollary~\ref{cor:tp_under_smallness} states that the turnpike occurs for small initial and terminal states under linear stabilizability and detectability. Relaxing the smallness conditions is one of major challenges in ($\mathrm{OCP}_2$). In \S~\ref{sctn:ex_2nd_class}, we will give a class of nonlinear systems for which the turnpike occurs for all initial states. This is done with the aid of the result in \cite{Sakamoto:20:prep} giving an estimates on the region for stable manifold in terms of nonlinear stabilizability. In the example in \S~\ref{sctn:ex_2nd_class}, the unstable manifold is linear and a geometric condition in Theorem~\ref{thm:trunpike_control_type2} is readily verified. 
\end{enumerate}
}
\end{rem}
\section{Examples}\label{sctn:examples}
\subsection{Problem $(\mathrm{OCP}_1)$}\label{sctn:ex_1st_class}
In this subsection, we show a class of nonlinear systems where the turnpike occurs in $(\mathrm{OCP}_1)$ for all target $z$. 
Let us consider the following class of nonlinear control systems
\begin{equation}
\begin{cases}
\dot x_1 = A_1x_1 +A_2(x_1,x_2)x_1\\
\dot x_2 = A_3x_2+B_2u,
\end{cases}\label{eqn:control_sys_ex1}
\end{equation}
where $A_1$ is an $n_1\times n_1$ Hurwitz matrix, $A_2:\mathbb{R}^{n_1}\times\mathbb{R}^{n_2}\to\mathbb{R}^{n_1\times n_1}$ is a $C^2$ function and $u\in\mathbb{R}^{m}$ is a control input. Assume that $(A_3,B_2)$ is stabilizable and $A_2(0,x_2)=0$ for all $x_2\in\mathbb{R}^{n_2}$. The cost function is 
\begin{equation}
J_1 = \frac{1}{2}\int_0^T |u|^2+|C_1x_1-z_1|^2+|C_2x_2-z_2|^2\,dt,\label{eqn:cost_f_ex1}
\end{equation}
where $C_1$, $C_2$ are constant matrices with appropriate dimensions, $(C_2,A_3)$ is detectable and $z_1\in\mathbb{R}^{r_1}$, $z_2\in\mathbb{R}^{r_2}$ are given constant vectors. 

The corresponding Hamiltonian system for this problem is 
\begin{equation}
\begin{cases}
\dot x_1= A_1x_1 +A_2(x_1,x_2)x_1\\ \dot x_2=A_3x_2-B_2B_2^\top p_2\\ \dot p_1= - C_1^\top(C_1x_1-z_1)-A_1^\top p_1 \\
\qquad \qquad -D_{x_1}\left[p_1^\top A_2(x_1,x_2)x_1\right]^\top \\
\dot p_2 =-C_2^\top(C_2x_2-z_2)\\
\qquad\qquad -A_3^\top p_2-D_{x_2}\left[p_1^\top  A_2(x_1,x_2)x_1\right]^\top.
\end{cases}\label{eqn:ham_trunpike_ex1}
\end{equation}
Using the stabilizability and detctability of $(C_3,A_3,B_2)$, it can be seen that there is an equilibrium $(0,x_{20}(z_2)$, $p_{10}(z_1),p_{20}(z_2))$ for (\ref{eqn:ham_trunpike_ex1}). At this equilibrium, the linearized matrix is
\[
\bmat{A_1&0&0&0\\0&A_3&0&-B_2B_2^\top\\
-C_1^\top C_1-\Gamma & 0&-A_1^\top&0\\ 0&-C_2^\top C_2&0&-A_3^\top},
\]
where $\Gamma=\Gamma(z_1,z_2)$ is an $n_1\times n_1$ symmetric matrix and therefore, it can be seen that it is a hyperbolic equilibrium. 

Let $P_1$, $P_3$ and $S_3$ be solutions for 
\begin{gather*}
P_1A_1 +A_1^\top +C_1^\top C_1+\Gamma=0,\\
P_3A_3+P_3A_3^\top -P_3B_2B_2^\top P_3 +C_2^\top C_2=0,\\
(A_3-B_2B_2^\top P_3)S_3+S_3(A_3-B_2B_2^\top P_3)^\top=B_2B_2^\top,
\end{gather*}
with $P_1=P_1^\top$, $P_3\geqslant0$, $S_3\geqslant0$ and $A_3-B_2B_2^\top P_3$ being Hurwitz. 
Using a linear coordinate transformation (see Appendix~\ref{sctn:appx}) 
\begin{gather*}
    \bmat{x_1\\x_2-x_{20}\\p_1-p_{10}\\p_2-p_{20}}=T\bmat{p_1'\\p_2'\\x_1'\\x_2'};\quad 
T=\bmat{I& 0& 0& 0\\ 0& I &0 &S_3\\ P_1& 0& I& 0\\0& P_3& 0& I+P_3S_3},\\
T^{-1}=\bmat{I& 0& 0& 0\\0& I+S_3P_3& 0& -S_3\\-P_1& 0& I& 0\\0& -P_3& 0& I},
\end{gather*}
the Hamiltonian system (\ref{eqn:ham_trunpike_ex1}) is rewritten as
\[
\begin{cases}
\dot{p}_1' = A_1p_1' +\psi_1(p_1',p_2',x_1',x_2')\\
\dot{p}_2'=(A_3-B_2B_2^\top P_3)p_2'+\psi_2(p_1',p_2',x_1',x_2')\\
\dot{x}_1'=-A_1^\top x_1' +\psi_3(p_1',p_2',x_1',x_2')\\
\dot{x}_2'=-(A_3-B_2B_2^\top P_3)^\top x_2'+\psi_4(p_1',p_2',x_1',x_2'),
\end{cases}
\]
where $\psi_j$, $j=1,\ldots,4$, are appropriately computed higher order terms. Since $\psi_j(0,0,x_1',x_2')$=0, $j=1,\ldots,4$, for all $x_1'$, $x_2'$, the unstable manifold $U$ at the equilibrium is the affine space $p_1'=p_2'=0$, or
\[
U=\{x_1=0, (I+S_3P_3)(x_2-x_{20})-S_3(p_2-p_{20})\}.
\]
Since $I+S_3P_3$ is nonsingular, which is shown using Lemma~\ref{lemma:PL+I} and Sylvester's determinant identity, for any $z_1$, $z_2$, $U$ intersects $p_1=p_2=0$ transversally. Now, using Theorem~\ref{thm:trunpike_control_type1}, for any $z_1$, $z_2$, if the initial point $(x_1(0),x_2(0))$ is close enough to $(0,x_{20})$, the optimal control for (\ref{eqn:control_sys_ex1})-(\ref{eqn:cost_f_ex1}) possesses the turnpike property. 

\begin{figure}[htp]
    \centering
    \includegraphics[keepaspectratio,width=0.45\textwidth]{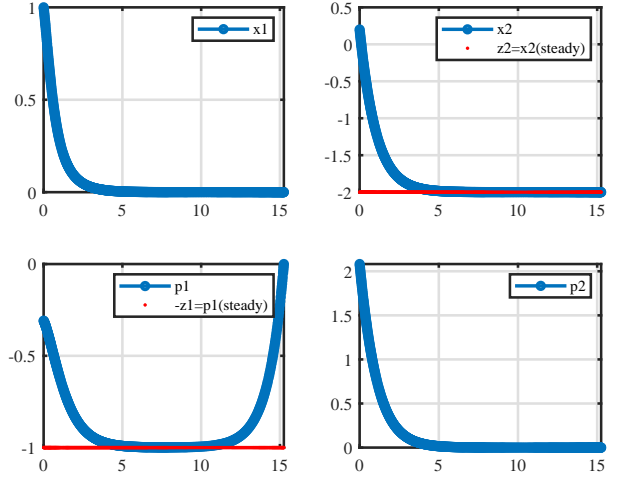}
    \caption{Optimal trajectory for (\ref{eqn:byrnes_ex}) $(x_1(0),x_2(0))=(1,0.2)$ and $(z_1,z_2)=(1,-2)$.}
    \label{fig:tpT2}
\end{figure}

As an example of the class of systems, a turnpike trajectory for a nonlinear optimal control problem 
\begin{subequations}\label{eqn:byrnes_ex}
\begin{gather}
    \dot{x}_1=-x_1+x_1^2x_2, \quad 
    \dot{x}_2=u\label{eqn:byrnes_ex_sys}\\
    J_1=\frac{1}{2}\int_0^Tu^2+(x_1-z_1)^2+(x_2-z_2)^2\,dt
\end{gather}
\end{subequations}
is depicted in Fig.~\ref{fig:tpT2}, where a solution of (SOP) is $(0,z_2,-z_1,0)$. %
\subsection{Problem $(\mathrm{OCP}_2)$}\label{sctn:ex_2nd_class}
Next, we show a class of nonlinear control systems for which estimates on (un)stable manifold of Hamiltonian systems obtained in \cite{Sakamoto:20:prep} are effective for the prediction of turnpike.  

Let us consider an $(n_1+n_2)$-dimensional system represented in Byrnes-Isidori normal form \cite{Byrnes:91:ieee-ac} for globally exponentially minimum phase nonlinear systems
\begin{equation}
\begin{cases}
    \dot{x} = q(x,y_1)\\
    \dot{y}_1=y_2\\
    \cdots\\
    \dot{y}_{n_2}=u,
\end{cases}\label{eqn:B-I_normalform}
\end{equation}
where $x\in\mathbb{R}^{n_1}$ and $q:\mathbb{R}^{n_1+1}\to\mathbb{R}^{n_1}$ is a smooth map with $q(0,0)=0$. We assume that $\dot x=q(x,0)$ is globally exponentially stable. It is known that (\ref{eqn:B-I_normalform}) is globally exponentially stabilizable via a smooth feedback. Therefore, representing $y=(y_1,\ldots,y_{n_2})$, for a cost functional 
\[
L(u,x,y)=\frac{1}{2}(|u|^2+|C_1x|^2+|C_2y|^2),
\]
the associated Hamiltonian system is hyperbolic at the origin if $C_2$ and the matrix defining $y$-dynamics is a detectable pair. If, in addition, $|C_1x|^2$ and $q(x,y_1)$ satisfy the growth condition in Proposition~\ref{prop:existence_OC_appdx}-(iv) with respect to $x$, the stable manifold $S$ of the Hamiltonian system satisfies $\pi_1(S)=\mathbb{R}^{n_1+n_2}$. Therefore, from Corollary~\ref{cor:tp_under_smallness}, the OCP has a solution for all $x_0$ and for sufficiently small $x_f$ that exhibits turnpike if the zero-state detectability condition is satisfied and $T$ is taken large enough. 

\begin{figure}[htp]
    \centering
    \includegraphics[keepaspectratio,width=0.45\textwidth]{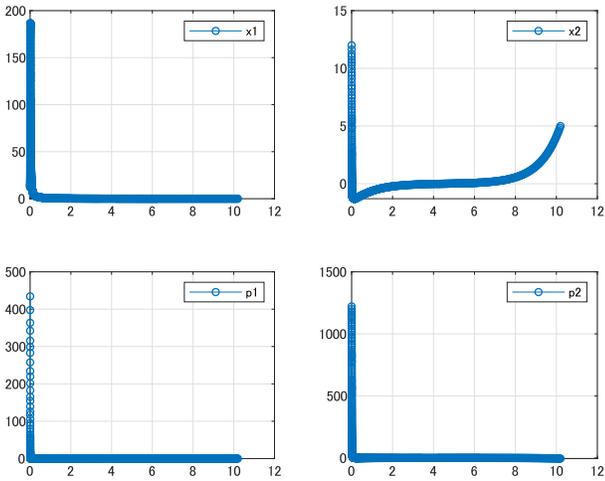}
    \caption{Optimal trajectory for (\ref{eqn:byrnes_ex_sys})-(\ref{eqn:cost_byrns_ex_type1}) with  $x_0=(12,12)$ and $x_f=(0,5)$.}
    \label{fig:tpT1_full}
\end{figure}

\begin{figure}[htp]
    \centering
    \includegraphics[keepaspectratio,width=0.45\textwidth]{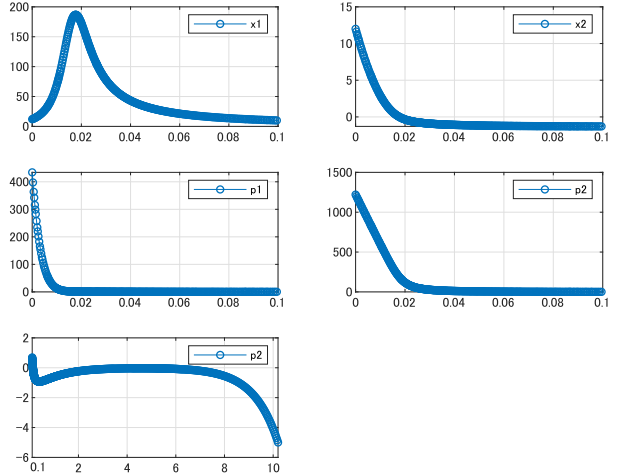}
    \caption{Enhanced figures of Fig.~\ref{fig:tpT1_full} for time spans $[0,0.1]$ and $[0.1,10]$~(last figure).}
    \label{fig:tpT1_eh}
\end{figure}

As a numerical example, consider (\ref{eqn:byrnes_ex_sys}), which is in Byrnes-Isidori normal form (see e.g., \cite{Byrnes:91:ieee-ac}), with 
\begin{equation}
J_2=\frac{1}{2}\int_0^T u^2+x_1^2+x_2^2\,dt.\label{eqn:cost_byrns_ex_type1}
\end{equation}
Introducing a cut-off function on $x_2$, the result in \cite{Sakamoto:20:prep} is applied to confirm that the turnpike occurs for all initial condition $x_0=(x_1(0),x_2(0))$ and terminal states $x_f\in \pi_1(U)$, where $U$ is the unstable manifold of the Hamiltonian system at the origin. Similarly to the previous subsection, $U$ is described as 
\[
U=\{ x_1=0,\ x_2+p_2=0 \}.
\]
Figs.~\ref{fig:tpT1_full}, \ref{fig:tpT1_eh} show the turnpike trajectory of the optimal control problem (\ref{eqn:byrnes_ex_sys})-(\ref{eqn:cost_byrns_ex_type1}) with $x_0=(12,12)$, $x_f=(0,5)$. In Fig.~\ref{fig:tpT1_eh}, $x_1(t)$, $x_2(t)$, $p_1(t)$, $p_2(t)$ are depicted for $t\in [0,0.1]$ while the last figure shows $p_2(t)$ for $t\in[0.1,10]$. From these figures, one sees that starting from $x_0=(12,12)$ at $t=0$, the states and costates rapidly grow during the time span $[0,0.02]$ and go to the steady optimal (the origin) by the time $t=0.1$ and then, the states reach the destination $x_f=(0,5)$ at $t=10$. The peak of this growth increases as $|x_0|$ increases. This growth of states is called "peaking phenomenon" of nonlinear stabilization \cite{Sussmann:91:ieee_tac} and it is interesting to see that peaking phenomenon appears in turnpike trajectory. 


\section{Discussions}\label{sctn:disscussion}
The geometric approach proposed in the present paper may be applied to more general cases where turnpike phenomena need more sophisticated analyses. Here, we discuss two kinds of extensions. 
\subsection{ Global analysis when (SOP) admits multiple solutions}
When (SOP) admits multiple solutions, multiple equilibria appear in associated Hamiltonian systems. If they are all hyperbolic, the $\lambda$-lemma still applies to draw pictures of flows around stable and unstable manifolds that are separatrices dividing the phase space (see, e.g., \cite[p.87 Corollary~1]{Palis:82:GTDS}). 

\begin{figure}[htp]
    \centering
    \includegraphics[keepaspectratio,width=0.47\textwidth]{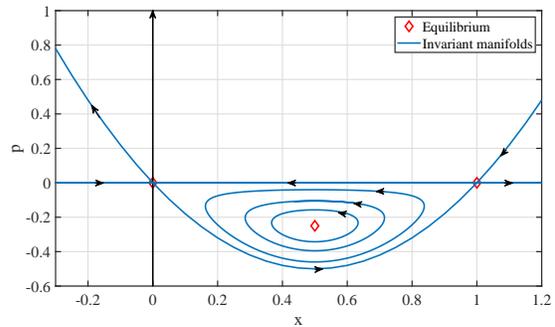}
    \caption{Invariant manifolds of Hamiltonian system for (\ref{eqn:multiple_equi}).}
    \label{fig:inv_manifolds_multiequili_EX}
\end{figure}

\begin{figure}[htp]
    \centering
    \includegraphics[keepaspectratio,width=0.47\textwidth]{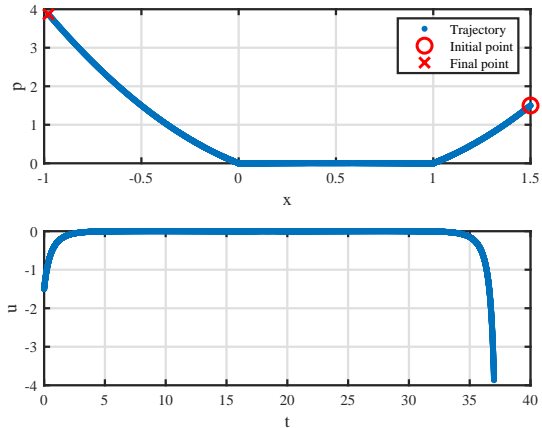}
    \caption{Prolonged turnpike due to two equilibrium points.}
    \label{fig:Prolonged_TP}
\end{figure}

\begin{figure}[htp]
    \centering
    \includegraphics[keepaspectratio,width=0.47\textwidth]{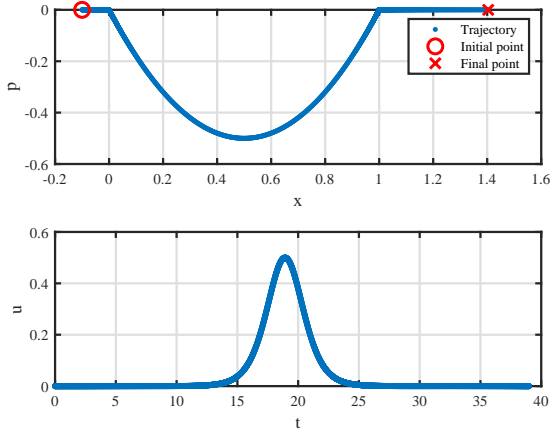}
    \caption{Bump turnpike due to two equilibrium points.}
    \label{fig:Bump_TP}
\end{figure}

Let us consider $(\mathrm{OCP}_2)_T$ for 
\begin{subequations}\label{eqn:multiple_equi}
\begin{gather}
    \dot x= -x+x^2+u,\\
    J_2(u)=\frac{1}{2}\int_0^Tu^2\,dt.
\end{gather}
\end{subequations}
Associated Hamiltonian system has three equilibrium points; $(0,0)$, $(1,0)$ and $(1/2,-1/4)$, the first two of which are the global solution of (SOP) and hyperbolic. Fig.~\ref{fig:inv_manifolds_multiequili_EX} shows stable and unstable manifolds, closed orbits around $(1/2,-1/4)$ and heteroclinic orbits connecting $(0,0)$ and $(0,1)$.  
From this figure and using the geometric method in the present paper, one immediately sees that for any initial point $x(0)$ and final point $x_f$, solution for $(\mathrm{OCP}_2)_T$ with large $T$ exists. For instance, trajectory in $x$-$p$ space and corresponding optimal input are depicted in Fig.~\ref{fig:Prolonged_TP} for $x(0)=1.5$, $x_f=-1$. Although the input response looks like turnpike, the response of $x$ for $u\sim0$ is not stationary but steady motion with nonzero velocity. 
Fig.~\ref{fig:Bump_TP} shows optimal trajectory and control response for $x(0)=-0.1$, $x_f=1.4$. Nonzero control is necessary to drive $x$ against stable vector field. The ratio of the time duration for nonzero control for the overall horizon can be arbitrarily small as $T\to\infty$ and in this sense, this can be also considered turnpike phenomenon. 

As for $(\mathrm{OCP}_1)$, when multiple global minimizers exist, 
an interesting question is raised in \cite{Phigin:20:nonunique_arxiv} as to which minimizer attracts turnpike for wider initial conditions. It is interesting to study how the geometry of these invariant manifolds affects turnpike occurrence and its strength in terms of the question. 
\subsection{Non-hyperbolic Hamiltonian systems} 
In \cite{Faulwasser:19:nolcos}, a concept of velocity turnpike or time-varying turnpike arising in mechanical systems is proposed combining trim primitives and turnpike properties. Motivated by that, the authors in \cite{Phigin:20:lackobs_arxiv} consider turnpike properties when detectability (ovservability) is not satisfied. A common feature in these cases is that associated Hamiltonian systems have zero eigenvalues. It is then interesting to consider the application of the $\lambda$-lemma for normally hyperbolic invariant manifolds \cite{Cresson:15:reg_chaotic} combining the classification result on Hamiltonian and symplectic matrices \cite{Laub:74:celes_mech}. 
\section{Conclusions}
In this paper, using techniques from dynamical system theory such as invariant manifolds and the $\lambda$-lemma, we showed that turnpike-like behavior naturally appears in hyperbolic dynamical systems. This is then applied to analyze Hamiltonian systems describing controlled trajectories to obtain sufficient conditions for optimal controls yielding the turnpike to exist. The framework proposed in the paper is geometric and an alternative to existing ones. Using the framework, we showed classes of nonlinear systems for which target $z$ or initial states can be taken arbitrarily large. 

Since our interests were to discover geometric nature in turnpike, we focused on OCPs without constraints and exponential turnpike. 
Future works include applications of this approach to more specific problems and considering OCPs with constraints, for which we mention an attempt to analyze turnpike in the maximum hands-off control \cite{sakamoto:20:cdc}. 
\vskip 1ex
{\em Acknowledgement.} The authors would like to thank Emmanuel Tr\'elat and Lars Gr{\"u}ne for their comments on the early version of the manuscript. The authors are also grateful to Dario Pighin for valuable discussions. 
\appendix

$\hfill${\bf Appendix}$\hfill$
\section{Results related with Riccati equations and linear Hamiltonian systems}\label{sctn:appx}
Let us consider Riccati equation 
\begin{equation}
PA+A^\top P-PRP+Q=0,\label{eqn:ric_appendix}
\end{equation}
where $A,R,Q\in\mathbb{R}^{n\times n}$ are constant matrices with $R,Q\geqslant0$. Suppose that $(A,R)$ is stabilizable and $(Q,A)$ is detectable. The following are known (see, e.g., \cite{Francis:86:CHCT,Lukes:69:sicon,Sakamoto:02:sicon}). 
\begin{lemma}\label{lemma:ham_appdx}
\begin{enumerate}[(i)]
    \item There is a solution $P\geqslant0$ to (\ref{eqn:ric_appendix}) such that $A_c:=A-RP$ is Hurwitz. 
    \item Let $L\leqslant0$ be a solution to a Lyapunov equation
    \[
    LA_c^\top +A_cL=R,
    \]
    then $\left[\begin{smallmatrix}I& L\\P&PL+I\end{smallmatrix}\right]$ is a symplectic matrix and its inverse is $\left[\begin{smallmatrix}LP+I& -L\\-P&I\end{smallmatrix}\right]$.
    \item 
    The Hamiltonian matrix $\mathrm{Ham}=\left[\begin{smallmatrix}A&-R \\-Q&-A^\top \end{smallmatrix}\right]$ is block-diagonalized as
    \begin{equation}
\mathrm{Ham}\bmat{I& L\\P&PL+I}=\bmat{I& L\\P&PL+I}\bmat{A_c& 0 \\0 &-{A_c}^\top }.\label{eqn:eigen_eq_appdx}
\end{equation}
\end{enumerate}
\end{lemma}
 The following lemma can be considered as dual version of Theorem~2 in \cite[p.90]{Francis:86:CHCT}, for which simplified proofs are given for the sake of self-containedness. 
\begin{lemma}\label{lemma:PL+I}
$PL+I$ is nonsingular. If, in addition, $(A,R)$ is controllable, then $L<0$ (negative definite). 
\end{lemma}
\quad{\em Proof.} Let $V:=PL+I$. From (\ref{eqn:eigen_eq_appdx}) we have 
\begin{subequations}\label{eqn:PL+I_proof}
\begin{align}
AL-RV &= -LA_c^\top\label{eqn:PL+I_proof1}\\
-QL-A^\top V &= -VA_c^\top\label{eqn:PL+I_proof2}
\end{align}
\end{subequations}
We show that the condition $\dim \mathrm{Ker}\,V\geqslant1$ leads to a contradiction. It can be shown from (\ref{eqn:PL+I_proof}) that $0\neq v\in\mathrm{ker}\,V$ satisfies $QLv=0$, $VA_c^\top v=0$ using $LV=V^\top L$ and $Q\geqslant0$, showing that $\mathrm{Ker}\,V$ is $A_c^\top$-invariant. Thus, we may assume that $v$ is an eigenvector of $A_c^\top$ with eigenvalue $\lambda$ with $\mathrm{Re}\,\lambda<0$. From (\ref{eqn:PL+I_proof2}), we have $ALv=-A_c^\top v=-\lambda Lv$ and therefore $(-\lambda I-A)Lv=0$. This shows that $\left[\begin{smallmatrix}Q\\-\lambda I-A\end{smallmatrix}\right]Lv=0$. With $\mathrm{Re}\,(-\lambda)>0$, the detectability of $(Q,A)$ implies $Lv=0$. This shows that $\left[\begin{smallmatrix}L\\V\end{smallmatrix}\right]v=0$ with $v\neq0$, which contradicts to Lemma~\ref{lemma:ham_appdx}(ii). 
The second statement can also be proved in a similar way, deriving $\lambda u=A^\top u$, $Ru=0$ for $0\neq u\in \mathrm{Ker}\,L$ and a contradiction. 
$\hfill\blacksquare$
\begin{lemma}\label{lemma:Phi_11_appdx}
Let 
\[
\bmat{\Phi_{11}(t)&\Phi_{12}(t)\\\Phi_{21}(t)&\Phi_{22}(t)}
=\exp[t\mathrm{Ham}],
\]
where $\Phi_{ij}(t)$, $i,j=1,2$, are $n\times n$ matrix functions of $t$. When $(A,R)$ is stabilizable and $(Q,A)$ is detectable, $\Phi_{11}(t)$ is nonsingular for $t\geqslant0$. 
\end{lemma}
\quad{\em Proof.} Using (\ref{eqn:eigen_eq_appdx}), 
\begin{align*}
\Phi_{11}(t) &= \exp[tA_c]\\
&\quad \times\left\{ I+\left( L-\exp[-tA_c]L\exp[-tA_c^\top] \right)P \right\}\\
&= \exp[tA_c](I+\tilde{L}(t)P),
\end{align*}
where we have set $\tilde{L}:=L-\exp[-tA_c]L\exp[-tA_c^\top]$. Since $\tilde{L}(0)=0$ and
\begin{align*}
\frac{d}{dt}\tilde{L}(t) =& \exp[-tA_c](A_cL+LA_c^\top)\exp[-tA_c^\top]\\
&=\exp[-tA_c]R\exp[-tA_c^\top]\geqslant0
\end{align*}
by Lemma~\ref{lemma:ham_appdx}(ii), $\tilde{L}(t)\geqslant0$ for $t\geqslant0$. If $\Phi_{11}(t)\eta=0$ for some $t\geqslant0$ and $\eta\in\mathbb{C}^n$, then we have $(I+\tilde{L}(t)P)\eta=0$ and therefore $\eta^\ast P\eta+\eta^\ast P\tilde{L}(t)P\eta=0$. This implies $P\eta=0$, $\tilde{L}(t)P\eta=0$ and we have $\eta=0$. $\hfill\blacksquare$
\section{Existence of infinite horizon optimal control and stable manifold of Hamiltonian systems}
This appendix introduces a result in \cite{Sakamoto:20:prep} on the existence of infinite horizon optimal control. The main result in the paper is under simpler growth conditions than those given below, but is more restrictive to apply. 

Let $U\subset\mathbb{R}^n$ be an open set containing the origin. 
A nonlinear system (\ref{eqn:nsys_general}) is said to be {\em $C^1$-exponentially stabilizable in $U$} if there exists a $C^1$ feedback control $u=k(x)$ with $k(0)=0$ such that the the closed loop system is exponentially stable with respect to $U$. Let $h(x)$ be a $C^1$ nonnegative function defined in $\mathbb{R}^n$. A system (\ref{eqn:nsys_general}) with output $y=h(x)$ is {\em zero-state detectable for $U$}, or simply {\em $(f,h)$ is zero-state detectable for $U$}, if the following holds. If a solution $x(t)$ with $x(0)\in U$ satisfies $h(x(t))=0$ for $t\geqslant 0$, then $x(t)\to0$ as $t\to\infty$. 

For system (\ref{eqn:nsys_general}), let $x=(x_1,x_2)$ with $x_1\in\mathbb{R}^{n_1}$, $x_2\in\mathbb{R}^{n_2}$, $n_1+n_2=n$ and rewrite it as
\begin{align*}
\frac{d}{dt}\bmat{x_1\\x_2} &=f(x_1,x_2)+g(x_1,x_2)u\\
    &=\bmat{f_1(x_1,x_2)\\f_2(x_1,x_2)}+\bmat{g_1(x_1,x_2)\\g_2(x_1,x_2)}u,
\end{align*}
where $f_j:\mathbb{R}^n\to\mathbb{R}^{n_j}$, $g_j:\mathbb{R}^n\to\mathbb{R}^{n_j\times m}$, $j=1,2$. 
Let $\varphi_R:\mathbb{R}^{n_2}\to\mathbb{R}$ be a $C^\infty$ cutoff function such that $\varphi_R(x_2)=1$ for $|x_2|<R$ and $\varphi_R(x_2)=0$ for $|x_2|\geqslant R+1$. Define $\tilde{f}_R(x_1,x_2):=f(x_1,\varphi_R(x_2)x_2)$ and $\tilde{g}_R(x_1,x_2):=g(x_1,\varphi_R(x_2)x_2)$. 

\begin{assumption}\label{assm:appendix}
\begin{enumerate}[(i)]
    \item System (\ref{eqn:nsys_general}) is $C^1$-exponentially stabilizable in $\Omega$, where $\Omega$ is an open set in $\mathbb{R}^n$ containing the origin.
    \item For a nonnegative $C^1$ function $h(x)$, there exist positive constants $p$, $\rho$, $c_h$ such that $h(x)\geqslant c_h|x|^p$ for $|x|>\rho$.  
    \item The pair $(f,h)$ is zero-state detectable for an open set containing $|x|\leqslant\rho$.
    \item For any $R>0$, there exist constants $c_f>0$, $c_g>0$, $0\leqslant\theta<1$, which may depend on $R$, such that 
\begin{align*}
    |\tilde{f}_R(x)|&\leqslant c_f|x|^{p+\theta}, \\
\|\tilde{g}_R(x)\|&\leqslant c_g|x|^{p/2+\theta},
\end{align*}
for sufficiently large $x\in\mathbb{R}^n$.
\item There exist constants $c_{f2}>0$, $c_{g2}>0$ and $0\leqslant\theta_2<1$ such that
\begin{align*}
    |{f}_2(x_1,x_2)|&\leqslant c_{f2}|x_2|^{p+\theta_2}, \\
\|{g}_2(x_1,x_2)\|&\leqslant c_{g2}|x_2|^{p/2+\theta_2},
\end{align*}
for all $x_1\in\mathbb{R}^{n_1}$ and sufficiently large $x_2\in\mathbb{R}^{n_2}$.
\end{enumerate}
\end{assumption}
\begin{proposition}\label{prop:existence_OC_appdx}
Under Assumption~\ref{assm:appendix}, for OPC (\ref{eqn:nsys_general}) and  
\begin{equation}
J=\int_0^\infty |u(t)|^2 +h(x(t))\,dt,\label{cost:appendix}
\end{equation}
there exists an optimal control for $x(0)\in\Omega$. Furthermore, for a Hamiltonian system associated with OCP (\ref{eqn:nsys_general})-(\ref{cost:appendix}), a stable manifold $S$ at the origin exists with the projection property $\Omega\subset\pi_1(S)$.
\end{proposition}

\bibliographystyle{plain} 
\bibliography{biblio/OptimalControl.bib,biblio/Main.bib,biblio/mechatro_satellite_aircraft.bib}

\end{document}